\documentclass[11pt]{amsart}
\catcode`\@=11

\def\newmathop#1{\expandafter\gdef\csname #1\endcsname{\mathop{\rm #1}\nolimits}}
\def\newmathopl#1{\expandafter\gdef\csname #1\endcsname{\mathop{\rm #1}\limits}}

\def\defsubscript#1{\expandafter\gdef\csname i#1\endcsname{_{\scriptscriptstyle #1}}}
\count0=26\loop\expandafter\expandafter\expandafter\defsubscript\@Alph{\the\count0}
\advance\count0 by -1\ifnum\count0>1\repeat

\def\defcalligraphic#1{\expandafter\gdef\csname c#1\endcsname
{{\ifmmode\mathcal#1\else$\mathcal#1$\fi}\spaceifletter}}
\count0=26\loop\expandafter\expandafter\expandafter\defcalligraphic\@Alph{\the\count0}
\advance\count0 by -1\ifnum\count0>1\repeat

\def\beq{\futurelet\comingchar\dobeq}     
  \def\dobeq{\ifx\comingchar*\begin{eqnarray*}\let\skipstar=\skiptoken
             \else\begin{eqnarray}\let\skipstar=\relax\fi\skipstar}
\def\eeq{\futurelet\comingchar\doeeq}     
  \def\doeeq{\ifx\comingchar*\end{eqnarray*}\let\skipstar=\skiptoken
             \else\end{eqnarray}\let\skipstar=\relax\fi\skipstar}
\def\skiptoken#1{}

\def\mdef#1#2{\def#1{\relax\ifmmode{#2}\else$#2$\fi\spaceifletter}}
\def\tdef#1#2{\def#1{\relax\ifmmode{\mathop{\rm #2}}\else#2\fi\spaceifletter}}

\def\spaceifletter{\futurelet\comingchar\dospaceifletter}
\def\dospaceifletter{\relax\ifmmode\else
  \ifcat A\noexpand\comingchar{} \fi
  \ifcat 0\noexpand\comingchar
  \ifx 0\noexpand\comingchar{} \fi
  \ifx 1\noexpand\comingchar{} \fi\ifx 2\noexpand\comingchar{} \fi
  \ifx 3\noexpand\comingchar{} \fi\ifx 4\noexpand\comingchar{} \fi
  \ifx 5\noexpand\comingchar{} \fi\ifx 6\noexpand\comingchar{} \fi
  \ifx 7\noexpand\comingchar{} \fi\ifx 8\noexpand\comingchar{} \fi
  \ifx 9\noexpand\comingchar{} \fi\fi
  \ifcat $\noexpand\comingchar{} \fi
  \ifcat \noexpand\relax\noexpand\comingchar{} \fi
\fi}

\def\isinlist#1#2#3{
\def\xisinlist##1##2#2##3##4\skiptoken{\ifx##3\skiptoken
  \def\isinlistfound{F}\else\def\isinlistfound{T}##1\fi
  }\xisinlist{#3}#1\donothing#2\skiptoken\skiptoken\donothing}

\def\inputdocument#1{\begingroup
  \def\documentstyle{\let\oldinput=\input\let\input=\skiptoken
  \let\inputifexists=\skiptoken
  \futurelet\comingchar\andocstyle}%
  \let\documentclass\documentstyle
  \def\andocstyle{\ifx\comingchar[\let\skipdoc=\skipdocstyle\else\let
  \skipdoc=\skiptoken\fi\skipdoc}%
  \def\document{\let\input=\oldinput\let\inputifexists=\oldinputifexists}%
  \def\enddocument{}%
  \input{#1}\endgroup}
\def\skipdocstyle[#1]#2{}

\catcode`\@=12

\mdef\P{\mathbb P}
\mdef\Q{\mathbb Q}
\mdef\R{\mathbb R}
\mdef\Z{\mathbb Z}
\mdef\C{\mathbb C}
\mdef\F{\mathbb F}

\let\tensor\otimes

\mdef\KF{K_2(\Q(C))}                           
\mdef\KC{K_2^T(C)}                             
\mdef\KCtors{K_2^T(C)/{\rm torsion}}           
\mdef\KCQ{K_2^T(C)\tensor\Q}                   
\mdef\KZ{K_2^{T,\Z}(C)}                        

\newmathop{Jac}
\newmathop{div}
\newmathop{Div}
\newmathop{Spec}
\newmathop{rk}
\newmathop{Pic}
\newmathop{lcm}

\def\Definition#1\par{\begin{definition}#1\end{definition}\par}
\def\Conjecture#1\par{\begin{conjecture}#1\end{conjecture}\par}
\def\Construction#1\par{\begin{construction}#1\end{construction}\par}
\def\Example#1\par{\begin{example}#1\end{example}\par}
\def\Proposition#1\par{\begin{proposition}#1\end{proposition}\par}
\def\Theorem#1\par{\begin{theorem}#1\end{theorem}\par}
\def\Remark#1\par{\begin{remark}#1\end{remark}\par}
\def\Lemma#1\par{\begin{lemma}#1\end{lemma}\par}
\def\Notation#1\par{\begin{notation}#1\end{notation}\par}

\def\ben{\begin{enumerate}}
\def\een{\end{enumerate}}
\def\bit{\begin{itemize}}
\def\eit{\end{itemize}}

{\catcode`\@=11\global\let\c@equation=\c@theorem}

\def\donothing{}         
\def\skiptoken#1{}

\def\prcounter#1{\ifx\donothing#1\donothing
  \gdef\@lastcounter{}%
  \gdef\prrfcmd{}%
\else
  \gdef\@lastcounter{\csname the#1\endcsname}%
  \gdef\prrfcmd{\refstepcounter{#1}}\addtocounter{#1}{1}
  \csname the#1\endcsname\addtocounter{#1}{-1}%
\fi}

\def\prindeedlabel#1]#2{{\label{#1}}\gdef\@lastproclaimlabel{#1}#2}

\def\prmaybelabel#1#2\par{\ifx[#1\prindeedlabel#2\par\else\gdef\@lastproclaimlabel{}#1#2\par\fi}

\def\genproclaim#1#2#3#4\par{{%
  \gdef\@lastproclaim{#1}%
  \medbreak\noindent{\bf#1\prcounter{#3}.\enspace}#2\prrfcmd\prmaybelabel#4\donothing\par
  \ifdim\lastskip<\medskipamount \removelastskip\penalty55\medskip\fi
}}

\catcode`\@=12
\newmathop{rad}
\def\d#1{#1}

\begin{document}

\title{LLL \& ABC}
\author{Tim Dokchitser}
\address{Dept. of Math. Sciences, University of Durham, UK}
\email{tim.dokchitser@durham.ac.uk}
\newcounter{equation}

\begin{abstract}
This note is an observation that the LLL algorithm applied to prime
powers can be used to find ``good'' examples for the ABC and Szpiro
conjectures.
\end{abstract}

\maketitle

Given non-zero integers $A$, $B$ and $C$, define
the {\em radical\/} $\,\rad(A,B,C)$ to be the product of primes dividing $ABC$,
the {\em size\/} to be $\,\max\,(|A|,|B|,|C|)\,$
and the {\em power\/} of the triple to be
$$
  P = P(A,B,C) = \frac{\log\max\,(|A|,|B|,|C|)}{\log\rad\,(A,B,C)} \>.
$$
The ABC conjecture of Masser-Oesterl\'e
states that for any real $\eta>1$, there are only finitely many
triples with $A,B,C$ relatively prime and $P(A,B,C)\ge\eta$ which satisfy
the ABC-equation
$$
  A + B = C \>.
$$
In particular,
some work has been done to look for examples of such solutions,
called ABC-triples, of as large power as possible (e.g. \cite{BB},
\cite{nit2}, \cite{weger}; see also the ABC page \cite{nitaj}
for an extensive list of references).
Similarly, one also looks for solutions with large Szpiro quotient
$$
  \rho = \rho(A,B,C) = \frac{\log\,|ABC|}{\log\rad\,(A,B,C)}\>.
$$
These relate via Frey curves to Szpiro's conjecture on the conductor
and the discriminant of elliptic curves. One can also formulate
the ABC conjecture over an arbitrary number field $K/\Q$
in terms of the generalized power \cite{broberg},
$$
  P(A,B,C) = \Bigl(\log\prod_\sigma\max\,(|A|_\sigma,|B|_\sigma,|C|_\sigma)\Bigr)
     \Bigm/
  \log\,\Bigl|\Delta_{K/\Q}\prod\limits_{p|ABC}N_{K/\Q}(p)\Bigr| \>.
$$
The product in the numerator is taken over all embeddings of $K$ into $\C$
and the denominator involves the discriminant of $K$ and the absolute norms
of the prime ideals dividing $A$, $B$ and~$C$.

The best examples known to date are: the ABC-triple found by E. Reyssat,
$$
  A = 2, \quad B = 3^{10}\!\cdot\!109, \quad C=23^5 \quad\text{with}\quad P=1.629...\>,
$$
and the Szpiro triple found by A. Nitaj,
$$
  A = 13\!\cdot\!19^6, \quad B = 2^{30}5, \quad C=3^{13}11^231 \quad\text{with}\quad \rho=4.419... \>.
$$
The best algebraic ABC example (by B. de Weger) involves
the roots $r_i$ of $x^3-2x^2+4x-4=0$,
$$
  (r_3\!-\!r_2)\,r_1^{52} + (r_1\!-\!r_3)\,r_2^{52} + (r_2\!-\!r_1)\,r_3^{52} = 0, \qquad P=1.920... \>.
$$
Following terminology of \cite{nitaj}, we call a triple with
$P>1.4$ a good ABC triple and one with $\rho>4$ a good Szpiro triple.

Our observation is that a simple method to look for such examples is
as follows.

Take, for instance, large prime powers $A_0\!=\!p^a$, $B_0\!=\!q^b$
and $C_0\!=\!r^c$ of comparable size with small $p,q$ and $r$.
Using the LLL lattice reduction algorithm \cite{LLL}, one can find
the smallest integral relation between them with respect to the $l^2$-norm,
\begin{equation}\label{rel}
  \alpha A_0 + \beta B_0 + \gamma C_0 = 0\>.
\end{equation}
Since the coefficients $\alpha, \beta$ and $\gamma$ are relatively small,
the numbers $\alpha p^a$, $\beta q^b$ and $\gamma r^c$ with suitably
chosen signs give a potential candidate for a high-powered ABC-triple.
Instead of prime powers one may also take a product of two prime powers or,
more generally, any large number with small radical. For example, the
smallest relation of the form (\ref{rel}) between
$$
  A_0 = 71^8, \quad B_0=2^55^{18}17^3\quad\text{and}\quad C_0 = 3^{38}
$$
has $(\alpha,\beta,\gamma)\!=\!(12649337, 336633577, -149459713)$.
This gives a previously unknown good triple,
$$
A = 71^8\,233^3, \quad B= 2^5\,5^{18}\,7^3\,17^3\,981439,\quad
    C=3^{38}\,13^4\,5233, \quad P=1.414... \>.
$$
Incidentally, this is the largest (with respect to size as above)
good ABC example known to the author.

To give an empirical analysis of this approach,
take $A_0\!=\!p^a$, $B_0\!=\!q^b$ and $C_0\!=\!r^c$ as above,
all three approximately of size $N$. Then in the worst case the smallest zero
combination of the form (\ref{rel}) has coefficients
$\alpha, \beta$ and $\gamma$ roughly of size
$\sqrt{3N}$. In fact, there are about $(\sqrt{3N})^3$ combinations
\begin{equation}\label{comb}
  i A_0 + j B_0 + k C_0, \qquad 0\le i,j,k\le\sqrt{3N}\>.
\end{equation}
As they are all of size at most $3\times N\sqrt{3N}=(\sqrt{3N})^3$, two of them must
be equal by the box principle and their difference gives a required relation.
Therefore, in the worst case the resulting triple has
$$
  P(A,B,C) \approx\frac{\log(N\sqrt{3N})}{3\log\sqrt{3N}+\log p+\log q+\log r} \>.
$$
For fixed $p, q$ and $r$, this expression tends to 1 as $N$ goes to infinity,
so the triples are ``on the edge'' of what is predicted by the ABC
conjecture.

It is easy to implement the above method to actually search
for some explicit new examples. To illustrate one practical consideration,
take
$$
  A_0 = 1, \quad B_0=3^4, \quad C_0=5^4\>.
$$
Here LLL reduction shows that the lattice of relations between these numbers
is generated by
$$
  v_1 = (23,-8,1) \quad\text{and}\quad v_2 = (12,23,-3)\>.
$$
Then $v_1$ gives a reasonable ABC-triple, but a few other small combinations
of $v_1$ and $v_2$ yield even better ones,
$$
  \begin{array}{lclll}
  v_1 &       = (23,-8,1) & \>\>\implies\>\> & (A,B,C) = (23, 5^4, 2^33^4), & \>\>\>P=0.990... \cr
  -v_1\!+\!2v_2 & = (1,54,-7) & \>\>\implies\>\> & (A,B,C) = (1, 2\!\cdot\!3^7, 5^47), & \>\>\>P=1.567... \cr
  4v_1\!-\!v_2 & = (104,-9,1) & \>\>\implies\>\> & (A,B,C) = (2^313, 5^4, 3^6), & \>\>\>P=1.104...\>. \cr
  \end{array}
$$
This suggests to run a search as follows. Take a list $L$ of numbers to serve
as $A_0, B_0$ and $C_0$. For example, choose bounds $M$ and $N$ and consider
all numbers less than $M$ whose prime factors are less than $N$. Then for
all distinct $A_0, B_0$, and $C_0$ in $L$ use LLL to determine the reduced
set of generators for the lattice of relations between them. Then take
``small'' combinations of these generators and check whether the resulting
ABC-triple has a sufficiently high power.

This does raise a question of how to decide which combinations
of the $v_i$ one should try.
A few experiments suggest that a sensible choice
is to try those $v\!=\!c_1v_1\!+\!c_2v_2$ which simply make one of the
elements of $v$ small.
This is a cheap way to keep the radical of the product of the
elements of $v$ to be as small as possible.
To achieve this for an index
$i\in\{1,2,3\}$, look at minus the $i$-th entry of $v_1$ divided by that
of $v_2$ and try several continued fraction approximations $c_2/c_1$
of this quotient. (In the example above, $-2/1$ is the first approximation
to $-23/12$ and it makes the first entry small.)
Incidentally, this is essentially using a form of LLL again.

Such a search has been implemented as a straightforward
Pari script with a running time of one weekend distributed over 30 PCs.
We found 145 out of 154 known good ABC-triples, 44 out of 47 known
good Szpiro triples and many new examples, listed in the two tables below.

Let us conclude with a few remarks.

First, note that although the approach presented here is apparently new,
LLL has been used in different ways in relation to the ABC conjecture;
see e.g. \cite{weger}.

Second, the method seems to works best when $A_0$, $B_0$ and $C_0$ are
of approximately equal size. In particular, one might expect it to be
more suitable for finding new Szpiro examples rather than ABC examples.
And, indeed, we have 48 new (47 known) good Szpiro triples but only
41 new (154 known) good ABC triples.

One can also perform a similar search in number fields and find
several interesting algebraic examples. For instance,
take the field $K\!=\!\Q(\sqrt{13})$ and $w=(3+\sqrt{13})/2$,
the fundamental unit of $K$. Let
$$
\begin{array}{lcr}
  A=w^{-5} (w\!-\!1)   & = &   0.0058594420567...\cr
  B=w^5 (w\!-\!2)      & = & 511.9941405579432...\cr
  C=2^9                & = & 512.0000000000000...
\end{array}
$$
Then $A$ and $B$ are conjugate and $A\!+\!B\!=\!C$. Moreover,
2 is prime in $K$ and $(w\!-\!1)(w\!-\!2)\!=\!3$. Hence the ABC-ratio is
$$
  P(A,B,C) = \frac{2\log(512)}{\log 13+\log3+\log3+\log4} = 2.0292288501126...
$$
So, this example is a new record for the algebraic ABC conjecture.

Finally, the same method can be also used to look for examples for the
generalization of the ABC-conjecture to solutions of $a_1+\ldots+a_n=0$,
known as the $n$-conjecture, see Browkin-Brzezinski \cite{BB}.

The author would like to thank N. Afanassieva for providing inspiration
to work on this problem and V. Dokchitser for helpful discussions.

\def\arraystretch{0.9}
\def\d#1{{\hbox{\fontsize{10pt}{10pt}\selectfont$#1$}}}

\begingroup
\begin{table}
$$
\begin{array}{|c|c|c|c|c|}
\hline
\d{A} & \d{B} & \d{C} & \d{\log_{10}\!C} & \d{\rho(A,B,C)} \cr
\hline
\vphantom{2^{X^{X^a}}}
\d{19^8\,43^4\,149^2} & \d{2^{15}\,5^{23}\,101} & \d{3^{13}\,13\!\cdot\!29^2\,37^6\,911} & \d{22.6} & \>\d{4.23181492}\>\cr
\d{2^7\,5^4\,7^{22}} & \d{19^4\,37\!\cdot\!47^4\,53^6} & \d{3^{14}\,11\!\cdot\!13^9\,191\!\cdot\!7829} & \d{23.9} & \d{4.21019250}\cr
\d{2^{17}\,3^{19}\,11\!\cdot\!25867} & \d{7^{12}\,23^7} & \d{5\!\cdot\!37^{10}\,53\!\cdot\!71} & \d{20.0} & \d{4.14980287}\cr
\d{7^8\,13\!\cdot\!89^3} & \d{3^{13}\,5^3\,11^4\,1499} & \d{2\!\cdot\!19^{12}} & \d{15.6} & \d{4.13636237}\cr
\d{11^3\,31^5\,101\!\cdot\!479} & \d{107^8} & \d{2^{31}\,3^4\,5^6\,7} & \d{16.3} & \d{4.13000150}\cr
\d{13^{10}\,37^2} & \d{3^7\,19^5\,71^4\,223} & \d{2^{26}\,5^{12}\,1873} & \d{19.5} & \d{4.12465150}\cr
\d{2^{55}\,23} & \d{3^{13}\,7^9\,13\!\cdot\!79^2} & \d{11^4\,43^6\,65353} & \d{18.8} & \d{4.10906942}\cr
\d{11^8\,13^9\,53} & \d{2^4\,5^{16}\,17\!\cdot\!547\!\cdot\!6163} & \d{7^6\,19^{12}} & \d{20.4} & \d{4.10809327}\cr
\d{233^4\,439} & \d{2^{15}\,3^{19}} & \d{5^8\,17^5\,71} & \d{13.6} & \d{4.10589886}\cr
\d{2^{13}\,71^2\,337^3} & \d{7^{13}\,1117^2} & \d{3^{21}\,13^3\,73^2} & \d{17.1} & \d{4.10470805}\cr
\d{5^7\,23^7\,1493} & \d{31^8\,3907^2} & \d{2^{52}\,3^2\,331} & \d{19.1} & \d{4.10115990}\cr
\d{2^{13}\,5^{12}\,13^4\,29} & \d{7^{16}\,19\!\cdot\!7451} & \d{3^{20}\,11^4\,353^2} & \d{18.8} & \d{4.08362226}\cr
\d{11^{11}\,73^2\,991\!\cdot\!306083} & \d{2^2\,3\!\cdot\!5^{11}\,7^{15}\,19^2} & \d{13^{15}\,31^5} & \d{24.2} & \d{4.08299029}\cr
\d{2\!\cdot\!5^9\,11^4\,41^2\,53^3} & \d{3^9\,7^{16}\,37} & \d{23^{11}\,40423} & \d{19.6} & \d{4.08262163}\cr
\d{31\!\cdot\!59^6} & \d{2^{25}\,3^{11}} & \d{5^3\,11^7\,13\!\cdot\!229} & \d{12.9} & \d{4.07920132}\cr
\d{3^{13}\,13\!\cdot\!23^3\,97^2} & \d{2^{37}\,157^2} & \d{5^5\,31^7\,67} & \d{15.8} & \d{4.07456723}\cr
\d{3^2\,73^{10}} & \d{5^{25}\,17^2\,23} & \d{2^2\,7^{11}\,827^2\,373357} & \d{21.3} & \d{4.07337395}\cr
\d{3^4\,5^{18}\,71\!\cdot\!419\!\cdot\!876581} & \d{2^{17}\,13^2\,19^{15}} & \d{7^6\,11^{12}\,977^3} & \d{26.5} & \d{4.06888583}\cr
\d{2^{26}\,11^4\,7639} & \d{5^6\,23^{11}} & \d{3^{18}\,47^4\,7879} & \d{19.2} & \d{4.06704768}\cr
\d{5^{16}\,19^2} & \d{3^8\,7^3\,89^4} & \d{2^{28}\,11^2\,6043} & \d{14.3} & \d{4.06668234}\cr
\d{5^7\,7^7\,19^2\,107} & \d{2^{14}\,11^9\,97} & \d{3^{10}\,23^7\,31} & \d{15.8} & \d{4.06231271}\cr
\d{2^7\,3\!\cdot\!821^5} & \d{13^{16}} & \d{5^{12}\,101^2\,324697} & \d{17.9} & \d{4.06159736}\cr
\d{3^4\,23^6\,1013^2} & \d{2^{47}\,5^3\,19^2} & \d{7\!\cdot\!131^7\,1373} & \d{18.8} & \d{4.06076852}\cr
\d{3^5\,5^{16}\,19^3} & \d{7^5\,23^2\,233^4\,1321} & \d{2^{48}\,43^2\,67} & \d{19.5} & \d{4.06075266}\cr
\d{83^2\,107^6} & \d{3^{17}\,5^{10}\,23} & \d{2\!\cdot\!7^2\,13^9\,17^2\,131} & \d{16.6} & \d{4.05751107}\cr
\d{2^6\,23^{11}\,53\!\cdot\!121523} & \d{7\!\cdot\!11^{17}\,13^3\,89} & \d{3^6\,17\!\cdot\!311^8} & \d{24.0} & \d{4.05572551}\cr
\d{5^9\,2141^2} & \d{2^{17}\,11\!\cdot\!53^4} & \d{3^2\,7\!\cdot\!19^9} & \d{13.3} & \d{4.05435143}\cr
\d{2^9\,3^{12}\,17^9\,1049} & \d{5^{19}\,23\!\cdot\!83\!\cdot\!491^2\,761} & \d{7^{12}\,13^4\,19^8} & \d{24.8} & \d{4.05071167}\cr
\d{13^3\,17^2\,131^5} & \d{2^{24}\,7^4\,11\!\cdot\!29^3\,103} & \d{3^3\,5^{19}\,47^2} & \d{18.1} & \d{4.04634190}\cr
\d{5^{26}\,11^2\,19^2} & \d{2^{11}\,139\!\cdot\!401^6\,463} & \d{3^{28}\,7^2\,37^2\,67^2\,89} & \d{23.8} & \d{4.04303710}\cr
\d{2^{47}\,3^7\,13^2} & \d{19^{11}\,23\!\cdot\!67^2\,227} & \d{5\!\cdot\!11\!\cdot\!257^6\,419^2} & \d{21.4} & \d{4.04183984}\cr
\d{2^4\,5^3\,17^6\,19^2\,151} & \d{7^{11}\,257^3} & \d{31^9\,37^2} & \d{16.6} & \d{4.04112782}\cr
\d{3^9\,17^7\,67} & \d{5^7\,7^9\,19\!\cdot\!439} & \d{2^4\,13^{10}\,23^3} & \d{16.4} & \d{4.04071176}\cr
\d{5^2\,17^6\,61^2\,269} & \d{2^{12}\,11^{12}} & \d{3^{19}\,7^5\,13\!\cdot\!53} & \d{16.1} & \d{4.03689092}\cr
\d{5^{20}\,4021} & \d{2^{40}\,13\!\cdot\!17^3\,6763} & \d{3^6\,7^7\,11^3\,29^6} & \d{20.7} & \d{4.03545668}\cr
\d{2^{23}\,5^2\,17^8} & \d{3^2\,67\!\cdot\!743^6} & \d{13^9\,23^4\,34679} & \d{20.0} & \d{4.03484329}\cr
\d{23^5\,43^5\,397} & \d{3^{33}\,5^3} & \d{2^{27}\,7\!\cdot\!29\!\cdot\!73^3\,101} & \d{18.0} & \d{4.03482389}\cr
\d{3^7\,7^5\,17^2\,239441} & \d{2^9\,5^{13}\,11^4} & \d{61^9} & \d{16.1} & \d{4.03406000}\cr
\d{31^7\,113\!\cdot\!491^2} & \d{5^{13}\,11\!\cdot\!13\!\cdot\!19^8} & \d{2\!\cdot\!3^6\,7^{18}\,1249} & \d{21.5} & \d{4.03083567}\cr
\d{2^2\,5\!\cdot\!11^3\,23^4\,29^7} & \d{13^9\,71^4\,113^2} & \d{3^{34}\,214033} & \d{21.6} & \d{4.02756838}\cr
\d{2^7\,7\!\cdot\!139^5} & \d{11\!\cdot\!41\!\cdot\!131^6} & \d{3^{27}\,5\!\cdot\!61} & \d{15.4} & \d{4.02755513}\cr
\d{2^{28}\,101\!\cdot\!197^4} & \d{37^{11}\,653} & \d{5^{14}\,7^2\,11^2\,2083^2} & \d{20.2} & \d{4.02174827}\cr
\d{5^2\,13\!\cdot\!37^6\,13789} & \d{2^9\,7^8\,47^4} & \d{3^{21}\,19^5} & \d{16.4} & \d{4.02095059}\cr
\d{3\!\cdot\!5\!\cdot\!67^9} & \d{11^8\,13^3\,47^3\,73} & \d{2^{14}\,7^3\,41^6\,149} & \d{18.6} & \d{4.01929332}\cr
\d{2^{27}\,3\!\cdot\!13^2\,19^5} & \d{5^6\,31^6\,263^2} & \d{37^9\,8677} & \d{18.1} & \d{4.00826664}\cr
\d{71^8\,233^3} & \d{2^5\,5^{18}\,7^3\,17^3\,981439} & \d{3^{38}\,13^4\,5233} & \d{26.3} & \d{4.00747592}\cr
\d{5^{15}\,13^6\,23^2} & \d{2^{31}\,61\!\cdot\!271^2\,19157} & \d{3^{26}\,7^3\,67^3} & \d{20.4} & \d{4.00512378}\cr
\d{11^7\,41^4} & \d{5^2\,7^7\,13^4\,211} & \d{2^{15}\,3^{16}\,127} & \d{14.3} & \d{4.00133657}\cr
\hline
\end{array}
$$
\caption{New Szpiro examples with $\rho(A,B,C)\!>\!4$}
\label{table2}
\end{table}

\begin{table}
$$
\begin{array}{|c|c|c|c|c|}
\hline
\d{A} & \d{B} & \d{C} & \d{\log_{10}\!C} & \d{P(A,B,C)} \cr
\hline
\vphantom{2^{X^{X^a}}}
\d{13^{10}\,37^2} & \d{3^7\,19^5\,71^4\,223} & \d{2^{26}\,5^{12}\,1873} & \d{19.5} & \>\d{1.50943262}\>\cr
\d{19^8\,43^4\,149^2} & \d{2^{15}\,5^{23}\,101} & \d{3^{13}\,13\!\cdot\!29^2\,37^6\,911} & \d{22.6} & \d{1.44280331}\cr
\d{3\!\cdot\!5^6\,7^8\,53} & \d{167^9} & \d{2\!\cdot\!11^6\,193^4\,20551} & \d{20.0} & \d{1.43823826}\cr
\d{2^{26}\,11^4\,7639} & \d{5^6\,23^{11}} & \d{3^{18}\,47^4\,7879} & \d{19.2} & \d{1.43813867}\cr
\d{7^8\,13\!\cdot\!89^3} & \d{3^{13}\,5^3\,11^4\,1499} & \d{2\!\cdot\!19^{12}} & \d{15.6} & \d{1.43785988}\cr
\d{17^4\,19^6} & \d{41^{10}\,1559} & \d{2^{12}\,3^{15}\,5\!\cdot\!29\!\cdot\!1567^2} & \d{19.3} & \d{1.43654400}\cr
\d{3^4\,7^2\,41} & \d{2^{25}\,227^7} & \d{5^9\,11^8\,2489197589} & \d{24.0} & \d{1.43575084}\cr
\d{3^{17}\,809} & \d{2^{27}\,11^9} & \d{5\!\cdot\!7^4\,13^5\,59\!\cdot\!1097^2} & \d{17.5} & \d{1.43485873}\cr
\d{11\!\cdot\!103^8} & \d{2^{45}\,3^7\,29\!\cdot\!37\!\cdot\!1997} & \d{5^{11}\,7^{10}\,79\!\cdot\!389^2} & \d{23.2} & \d{1.43360120}\cr
\d{7^{11}\,19} & \d{5^{12}\,1019\!\cdot\!7151^2} & \d{2^{28}\,3^{12}\,11^3\,67} & \d{19.1} & \d{1.43309388}\cr
\d{2^{17}\,13^3} & \d{7^3\,11^7\,43^2\,5801} & \d{3^{17}\,17^6\,23} & \d{16.9} & \d{1.43234742}\cr
\d{7^2\,23^5} & \d{3^8\,11^{12}\,4703} & \d{2\!\cdot\!5^2\,13^2\,19^3\,29^6\,53^2} & \d{20.0} & \d{1.42991591}\cr
\d{2^{20}\,79\!\cdot\!97} & \d{5^3\,7^6\,11^{10}} & \d{3^4\,13^7\,8663^2} & \d{17.6} & \d{1.42942802}\cr
\d{29^4} & \d{2^{14}\,3^3\,31\!\cdot\!47^2\,199^3} & \d{7^{12}\,4153^2} & \d{17.4} & \d{1.42836105}\cr
\d{2^{27}\,809} & \d{5^7\,7^6\,13^5} & \d{3^{23}\,36251} & \d{15.5} & \d{1.42460741}\cr
\d{2^{13}\,3^{18}\,2069} & \d{13^3\,29^7\,271^3} & \d{5^{14}\,23^3\,3187^2} & \d{20.9} & \d{1.42340611}\cr
\d{31^7\,113\!\cdot\!491^2} & \d{5^{13}\,11\!\cdot\!13\!\cdot\!19^8} & \d{2\!\cdot\!3^6\,7^{18}\,1249} & \d{21.5} & \d{1.42308625}\cr
\d{3^4\,23^6\,1013^2} & \d{2^{47}\,5^3\,19^2} & \d{7\!\cdot\!131^7\,1373} & \d{18.8} & \d{1.42201537}\cr
\d{3\!\cdot\!5\!\cdot\!13^6} & \d{2^7\,7^5\,53^6\,2287} & \d{11^3\,37^7\,929^2} & \d{20.0} & \d{1.42137859}\cr
\d{233^4\,439} & \d{2^{15}\,3^{19}} & \d{5^8\,17^5\,71} & \d{13.6} & \d{1.42081322}\cr
\d{2^{13}\,71^2\,337^3} & \d{7^{13}\,1117^2} & \d{3^{21}\,13^3\,73^2} & \d{17.1} & \d{1.42075105}\cr
\d{17^2\,47^5\,73} & \d{2^{31}\,5^9} & \d{3\!\cdot\!13^7\,4723^2} & \d{15.6} & \d{1.41627640}\cr
\d{2^7\,5^4\,7^{22}} & \d{19^4\,37\!\cdot\!47^4\,53^6} & \d{3^{14}\,11\!\cdot\!13^9\,191\!\cdot\!7829} & \d{23.9} & \d{1.41582933}\cr
\d{5^{20}\,4021} & \d{2^{40}\,13\!\cdot\!17^3\,6763} & \d{3^6\,7^7\,11^3\,29^6} & \d{20.7} & \d{1.41575870}\cr
\d{3^{22}\,9787^2} & \d{5^{10}\,11\!\cdot\!29^{10}\,109} & \d{2^{37}\,89^3\,167^2\,1823} & \d{24.7} & \d{1.41570162}\cr
\d{71^8\,233^3} & \d{2^5\,5^{18}\,7^3\,17^3\,981439} & \d{3^{38}\,13^4\,5233} & \d{26.3} & \d{1.41457078}\cr
\d{29^4\,2213^2} & \d{3\!\cdot\!13^2\,23^{12}\,89\!\cdot\!14717} & \d{2^9\,5^{16}\,11^9\,79} & \d{25.2} & \d{1.41234761}\cr
\d{5^4\,17\!\cdot\!349^3} & \d{7^{17}\,109} & \d{2^{35}\,3^5\,3037} & \d{16.4} & \d{1.41227598}\cr
\d{3^{11}\,7^2\,37\!\cdot\!47} & \d{2^8\,17^3\,101^2\,191^4} & \d{5^2\,353^7} & \d{19.2} & \d{1.41143590}\cr
\d{2^{29}\,13} & \d{5^{11}\,269^3} & \d{3^5\,7^2\,17^6\,3307} & \d{15.0} & \d{1.41090947}\cr
\d{5^4\,53^2\,59^4\,101} & \d{2^4\,11\!\cdot\!23^{15}} & \d{3^{14}\,7^{12}\,463\!\cdot\!1531} & \d{22.7} & \d{1.41032306}\cr
\d{2^{11}\,3^4\,101^4\,29221} & \d{13^{19}} & \d{5^{15}\,17\!\cdot\!53093^2} & \d{21.2} & \d{1.40944818}\cr
\d{7^4} & \d{3^{21}\,11^2\,13^4\,138493} & \d{2^{26}\,5^7\,383\!\cdot\!1579^2} & \d{21.7} & \d{1.40900897}\cr
\d{11^3\,31^5\,101\!\cdot\!479} & \d{107^8} & \d{2^{31}\,3^4\,5^6\,7} & \d{16.3} & \d{1.40714951}\cr
\d{3^5\,5^{16}\,19^3} & \d{7^5\,23^2\,233^4\,1321} & \d{2^{48}\,43^2\,67} & \d{19.5} & \d{1.40487025}\cr
\d{5^7\,23^7\,1493} & \d{31^8\,3907^2} & \d{2^{52}\,3^2\,331} & \d{19.1} & \d{1.40479669}\cr
\d{2^2\,3^4\,163^3\,1006151} & \d{43^{13}} & \d{11^9\,29^4\,101^3} & \d{21.2} & \d{1.40308397}\cr
\d{2^{11}\,3^{17}\,13^2\,19^2} & \d{29\!\cdot\!41^2\,83^8} & \d{5^4\,47^2\,53\!\cdot\!107^6} & \d{20.0} & \d{1.40244119}\cr
\d{3^6\,7^4\,43\!\cdot\!16421} & \d{5^{12}\,439^6} & \d{2^{59}\,41\!\cdot\!73939} & \d{24.2} & \d{1.40168452}\cr
\d{7^9\,13^4} & \d{2^{10}\,23^3\,173^4} & \d{3^{12}\,5\!\cdot\!11^6\,2371} & \d{16.0} & \d{1.40127027}\cr
\d{17^4} & \d{2\!\cdot\!7^{12}\,29^3\,743} & \d{3^9\,5^6\,13^5\,23\!\cdot\!191} & \d{17.7} & \d{1.40004159}\cr

\hline
\end{array}
$$
\caption{New ABC examples with $P(A,B,C)\!>\!1.4$}
\label{table1}
\end{table}

\endgroup

\end{document}